\documentclass[11pt]{article}
\usepackage{amsmath,amsthm,amssymb,tikz}
\usepackage{hyperref}

\hypersetup{
pdfauthor = {Boris Bukh, Gabriel Nivasch},
pdfsubject = {Mathematics},
pdfkeywords = {epsilon-approximants, epsilon-nets, regularity lemma},
pdfstartview = {FitH},
pdfpagemode = {UseNone}}

\hypersetup{pdftitle = {One-sided epsilon-approximants}}

\title{One-sided epsilon-approximants}
\author{Boris Bukh\thanks{Department of Mathematical Sciences, Carnegie
Mellon University, Pittsburgh, PA 15213, USA. Research was supported in
part by Churchill College, Cambridge, by U.S. taxpayers through NSF grant
DMS-1301548, and by Alfred P.\ Sloan Foundation through Sloan Research
Fellowship.} \and Gabriel Nivasch\thanks{Department of Computer Science, Ariel University, Ariel, Israel.}}
\date{}

\theoremstyle{plain}
\newtheorem{theorem}{Theorem}
\newtheorem{lemma}[theorem]{Lemma}

\newtheorem{proposition}[theorem]{Proposition}
\newtheorem{corollary}[theorem]{Corollary}
                  
\newtheorem{claim}{Claim}                             

\newcommand*{\abs}[1]{\lvert #1\rvert}                
\newcommand*{\babs}[1]{\left\lvert #1\right\rvert}    
\newcommand*{\R}{\mathbb{R}}                          
\newcommand*{\N}{\mathbb{N}}                          
\newcommand*{\Hg}{\mathcal{H}}                         

\newcommand*{\eqdef}{\stackrel{\text{\tiny{def}}}{=}} 
\def\mygobble#1{}
\newcommand*{\veps}{\varepsilon}                      
\def\vepsapproximant/{$\veps$\nobreakdash-approx\-i\-mant}
\def\vepsapproximants/{$\veps$\nobreakdash-approx\-i\-mants}
\def\vepsapproximating/{$\veps$\nobreakdash-approx\-i\-mating}
\def\vepstwoapproximating/{$(\veps/2)$\nobreakdash-approx\-i\-mating}
\def\vepsnet/{$\veps$\nobreakdash-net}
\def\vepsnets/{$\veps$\nobreakdash-nets}

\newcommand*{\TODO}[1]{}                              
\newcommand*{\NB}[1]{}                                

\newcommand*{\odd}{\mathrm{odd}}                      
\newcommand*{\even}{\mathrm{even}}                    
\DeclareMathOperator{\sgn}{sgn}                       
\DeclareMathOperator{\conv}{conv}                     
\DeclareMathOperator{\ahull}{ahull}                   
\DeclareMathOperator{\orient}{orient}                 
\DeclareMathOperator{\VCdim}{VC-dim}                  
\DeclareMathOperator{\tw}{tw}                         
\DeclareMathOperator{\avg}{avg}                       
\DeclareMathOperator{\Tver}{Tver}                     
\DeclareMathOperator{\OT}{OT}                         

\newcommand*{\F}{\mathcal{F}}                         
\newcommand*{\Fconv}{\mathcal{F}_{\operatorname{conv}}}  

\begin{document}

\maketitle

\begin{center}
\textit{In memory of a great teacher,\\Jirka Matou\v{s}ek}
\end{center}

\begin{abstract}
  Given a finite point set $P\subset\R^d$, we call a multiset $A$ a \emph{one-sided weak $\varepsilon$-approximant} for $P$ (with respect
  to convex sets), if
  $\abs{P\cap C}/\abs{P}-\abs{A\cap C}/\abs{A}\leq\veps$ for every convex
  set~$C$.

  We show that, in contrast with the usual (two-sided) weak
  \vepsapproximants/, for every set $P\subset \R^d$
  there exists a one-sided weak \vepsapproximant/ of size bounded
  by a function of $\veps$ and~$d$.	
\end{abstract}

\section{Introduction}\label{sec_intro}

A common theme in mathematics is approximation of large, complicated objects by smaller, simpler
objects. This paper proposes a new notion of approximation in combinatorial geometry,
which we call one-sided \vepsapproximants/. It is a notion of approximation that is in strength
between \vepsapproximants/ and \vepsnets/. We recall these two notions first.

Let $P\subset \R^d$ be a finite set, and $\F\subset 2^{\R^d}$ a family of sets in $\R^d$.
In applications, the family $\F$ is usually a geometrically natural family, such as the family
of all halfspaces, the family of all simplices, or the family of all convex sets.
A finite set $A\subset \R^d$ is called an \emph{\vepsapproximant/ for $P$ with
respect to $\F$} if
\[
  \babs{\frac{\abs{C\cap P}}{\abs{P}}-\frac{\abs{C\cap A}}{\abs{A}}}\leq \veps\qquad\text{for all }C\in\F.
\]
The notion of an \vepsapproximant/ was introduced by Vapnik and Chervonenkis \cite{vapnik_chervonenkis}
in the context of statistical learning theory. They associated to each family $\F$ 
a number $\VCdim(\F)\in\{1,2,3,\dotsc,\infty\}$, which has become known as
\emph{VC dimension}, and proved that if $\VCdim(\F)<\infty$, then
every set $P$ admits an \vepsapproximant/ $A$ of size $\abs{A}\leq C_{\VCdim(\F)} \veps^{-2}$,
a bound which does not depend on the size of~$P$. The \vepsapproximants/
that they constructed had the additional property that $A\subset P$. Following tradition, we say that $A$ is a \emph{strong \vepsapproximant/} if $A\subset P$.
When we wish to emphasize that our \vepsapproximants/ are not necessarily subsets of $P$, we call them \emph{weak \vepsapproximants/}.
The bound has been improved 
to $\abs{A}\leq C_{\VCdim(\F)}\veps^{-2+2/(\VCdim(\F)+1)}$ (see \cite[Theorem~1.2]{matousek_vcdisc_upper} and \cite[Exercise 5.2.7]{matousek_discbook})
which is optimal \cite{alexander_vcdisc_lower}. 

In a geometric context, Haussler and Welzl \cite{haussler_welzl} introduced \vepsnets/. With
$P$ and $\F$ as above, a set $N$ is called an \emph{\vepsnet/ for $P$ with respect to $\F$}
if
\[
   \frac{\abs{C\cap P}}{\abs{P}}>\veps \implies  C\cap N\neq \emptyset\qquad\text{for all }C\in\F.
\]
An \vepsapproximant/ is an \vepsnet/, but not conversely. While an \vepsnet/ is a weaker
notion of approximation, its advantage over an \vepsapproximant/ is that every set~$P$
admits an \vepsnet/ of size only $C_{\VCdim(\F)} \veps^{-1}\log \veps^{-1}$, which is smaller
than the bound for the \vepsapproximants/. The \vepsnets/ constructed by Haussler and Welzl
are also strong, i.e., they satisfy $N\subset P$.

Most geometrically important families $\F$ have a bounded VC dimension. A notable exception
is the family $\Fconv$ of all convex sets. Indeed, it is easy to see that a
set of $n$ points in convex position does not admit any strong \vepsnet/ of size smaller than
$(1-\veps)n$ with respect to $\Fconv$. Alon, B\'ar\'any F\"uredi, and Kleitman \cite{alon_barany_furedi_kleitman}
showed that for every $P\subset \R^d$ there exists a (weak) \vepsnet/ of size bounded solely
by a function of $\veps$ and $d$. No extension of their result to \vepsapproximants/
is possible.
\begin{proposition}\label{prop_pessimism}
If $P\subset \R^2$ is a set of $n$ points in convex position, then every \vepsapproximant/
with respect to $\Fconv$ has size at least $n(\tfrac{1}{4}-\veps/2)$. 
\end{proposition}
\begin{proof}
Let $p_1,p_2,\dotsc,p_n$ be the enumeration of the vertices of $P$
in clockwise order along the convex hull of $P$. For $i=1,\dotsc,\lfloor (n-1)/2\rfloor$
write $T_i$ for the triangle $p_{2i-1},p_{2i},p_{2i+1}$. Suppose
$A\subset \R^2$ is an \vepsapproximant/ for~$P$. Let $I\eqdef\{i : T_i\cap A=\emptyset\}$.
Note that $\abs{I}\geq n/2-2\abs{A}-1$ since each point of~$A$ lies in at most two triangles. Define $S\eqdef \{p_1,p_3,p_5,\dotsc\}$ to be
the odd-numbered points, and let $S'\eqdef S\cup\{p_{2i} : i\in I\}$. Let
$C\eqdef\conv S$ and $C'\eqdef\conv S'$. Then $C\cap A=C'\cap A$,
but $\abs{C'\cap P}/\abs{P}-\abs{C\cap P}/\abs{P}=\abs{I}/\abs{P}>\veps$
if $\abs{A}<\abs{P}(\tfrac{1}{4}-\veps/2)$.
\end{proof}

In light of Proposition~\ref{prop_pessimism}, we introduce a new concept.
A multiset\footnote{In this paper we allow $A$ to be a multiset. While the results of
this paper continue to hold if we require $A$ to be a set, the proofs become more technical. We sketch
the necessary changes in the final section.} $A\subset\R^d$ is a \emph{one-sided \vepsapproximant/ for $P$ with respect to the family $\F$} if
\[
 \frac{\abs{C\cap P}}{\abs{P}}-\frac{\abs{C\cap A}}{\abs{A}}\leq \veps\qquad\text{for all }C\in\F.
\]
In other words, if $C\in\F$, then $C$ might contain many more points of $A$ than expected, but never much fewer.
It is clear that an \vepsapproximant/ is a one-sided \vepsapproximant/, and that 
a one-sided \vepsapproximant/ is an \vepsnet/.

Our main result shows that allowing one-sided errors is enough to 
sidestep the pessimistic Proposition~\ref{prop_pessimism}.

\begin{theorem}\label{thm_main}
Let $P\subset \R^d$ be a finite set, and let $\veps\in (0,1]$ be
a real number. Then $P$ admits a one-sided \vepsapproximant/ with respect to $\Fconv$ of size at most $g(\veps,d)$,
for some $g$ that depends only on $\veps$ and on~$d$.
\end{theorem}

Unfortunately, due to the use of a geometric Ramsey theorem, our bound on $g$ is very weak:
\[
  g(\veps,d)\leq \tw_d\bigl(\veps^{-c}\bigr)
\]
for some constant $c>1$ that depends only on $d$, where the tower function is given by $\tw_1(x)\eqdef x$ and $\tw_{i+1}(x)\eqdef 2^{\tw_i(x)}$. We believe
this bound to be very far from sharp.

In the rest of the paper we omit the words ``with respect to $\Fconv$'' when referring to
one-sided approximants.


\section{Outline of the construction and of the paper}
At a high level, the proof of Theorem~\ref{thm_main} can be broken into three steps:
\begin{enumerate}
\item We replace the given set $P$ by a bounded-size set $\hat{P}$. The price of this replacement
is an extra condition that a one-sided \vepsapproximant/ $A$ for $\hat{P}$ would need to satisfy
to be a one-sided \vepsapproximant/ for $P$. Namely, $A$ must be a one-sided \vepsapproximant/
for a \emph{semialgebraic reason}.

\item We break $\hat{P}$ into long \emph{orientation-homogeneous} subsequences $S_1,S_2,\dotsc,S_m$.

\item For each $S_i$ we give an explicit one-sided \vepsapproximant/ $A_i$ satisfying the semialgebraicity condition. 
The union of $A_1,\dotsc,A_m$ is then the desired \vepsapproximant/.
\end{enumerate}

Step~1 relies on the semialgebraic regularity lemma from \cite{fox_pach_suk_regularity}, which we recall in
Section~\ref{sec_regularity}. Given a fixed set
$\Phi$ of semialgebraic predicates, this lemma permits us to replace
$P$ with a constant-sized $\hat{P}$ that behaves similarly to $P$ with respect
to the predicates in $\Phi$. Since in steps~2 and~3 we employ only one predicate,
in our case we have $\abs{\Phi}=1$. We define that predicate in Section~\ref{sec_geometric}.

In step~1 we lose some control on the interaction between parts of $\hat{P}$. 
To remedy this, we use a well-known hypergraph Tur\'an theorem, discussed in Section~\ref{sec_hypergraph},
to extract well-behaved chunks from~$\hat{P}$.

The construction of $A_i$ in step~3 consists of Tverberg points of a certain family $\mathcal{F}$ of subsets of $S_i$.
The property that $\mathcal{F}$ needs to satisfy is most naturally described in terms of interval chains, which are
introduced in Section~\ref{sec_interval}.
The actual construction of requisite interval chains is based on the idea behind the regularity lemma for words from \cite{words_regularity,feige_repeat}.
To obtain a better quantitative bound, we eschew using the lemma directly and provide an alternative argument.
This is also done in Section~\ref{sec_interval}.

All the ingredients are put together in Section~\ref{sec_construction}. 

The paper concludes with several remarks and open problems.

\section{Geometric preliminaries}\label{sec_geometric}

The \emph{convex hull} of a point set $P$ is denoted $\conv P$,
and its \emph{affine hull} is denoted $\ahull P$.

Tverberg's theorem (see, e.g., \cite[p.~200]{matousek_geometry_book}) asserts that any set $Q\subset \R^d$ of $(s-1)(d+1)+1$ points
can be partitioned into $s$ pairwise disjoint subsets whose convex hulls intersect. We denote by $\Tver_s(Q)$ an arbitrary point in such an intersection. A special case of Tverberg's theorem is the case $s=2$, which is due to Radon \cite{radon_orig}. 
In that case, if $Q$ is in general position (no $d+1$ points are affinely dependent), then the partition is unique and
$\Tver_2(Q)$ is also unique. 

A \emph{(geometric) predicate} of arity $k$ is a property that a $k$-tuple of points $p_1, \ldots, p_k$ might or might not satisfy. 
A predicate is \emph{semialgebraic} if it is a Boolean combination of expressions of the form $f(p_1,\dotsc,p_k)\geq 0$, where the $f$'s are polynomials.
Predicates that depend on the sign of a single polynomial are especially useful, we call then \emph{polynomial predicates}. For brevity we will identify 
polynomial predicates with the underlying polynomials.

An important polynomial predicate is the \emph{orientation} of a $(d+1)$\nobreakdash-tuple of points in~$\R^d$. The orientation of $p_0,\dotsc,p_d\in \R^d$ is given by
\begin{equation*}
\orient(p_0,\dotsc,p_d) \eqdef \sgn \det
\begin{bmatrix}
p_0 & \cdots &p_d\\
1 & \cdots &1
\end{bmatrix}.
\end{equation*}
We have $\orient(p_0, \ldots, p_d) = 0$ if and only if the points are affinely dependent.

\subsection{Orientation-homogeneous sequences}\label{subsec:orientation_homog}

We will call a sequence of points in $\R^d$ \emph{orientation-homogeneous} if all its $(d+1)$-tuples have the same nonzero orientation.
It is well known that every orientation-homogeneous sequence is in convex position,
and that the convex hull of such a sequence is combinatorially equivalent to a cyclic polytope
(see, e.g., \cite{ziegler} for background).

Let $P=(p_1,\dotsc,p_n)$ be an orientation-homogeneous sequence. For a set $I=\{i_1<\dotsb<i_m\}$, define the subsequence of $P$ 
indexed by $I$ by $P_I\eqdef(p_{i_1},\dotsc,p_{i_m})$. If $\abs{I}=d$, then the $d$ points $P_I$ span
a hyperplane $H_I\eqdef \ahull P_I$ in $\R^d$. It is simple to tell to which side of $H_I$ a 
point $p_j\in P\setminus P_I$ belongs: The index set $I$ partitions $[n]\setminus I$ into $d+1$ 
intervals (some of which might be empty). The side of $H_I$ to which $p_j$ belongs
depends only on the parity of the interval number to which $p_j$ belongs. In other words, $p_j$ is on one side if
$j\in (-\infty,i_1)\cup (i_2,i_3)\cup\dotsb$, and on the other side
if $j\in (i_1,i_2)\cup(i_3,i_4)\cup\dotsb$. Hence, two points
$p_j$ and $p_{j'}$ with $j<j'$ lie on the same side of $H_I$ if and only if 
$[j,j']\cap I$ is of even size.

Of particular interest to us are sets $I$ of size $d+2$. We define, for such a set $I=\{i_1<i_2<\dotsb<i_{d+2}\}$,
a partition $I=I_\odd\cup I_\even$, where $I_\odd\eqdef\{i_1,i_3,\dotsc\}$ and $I_\even\eqdef\{i_2,i_4,\dotsc\}$. 
\begin{lemma}\label{lem_interleave_radon}
If $P$ is orientation-homogeneous and $\abs{I}=d+2$, then the convex sets $\conv P_{I_\odd}$ and $\conv P_{I_\even}$ intersect.
\end{lemma}
\begin{proof}
Indeed, suppose they are
disjoint, and hence there exists a hyperplane $H$ that separates $P_{I_\odd}$ from $P_{I_\even}$.
Then $H$ can be perturbed into a hyperplane $H'$ that goes through some $d$ points of $P_I$,
i.e., $H'=H_J$ for some $J\subset I$, $\abs{J}=d$. The set $P_{I\setminus J}$
consists of two points, say $p_i,p_{i'}$ with $i<i'$, and they belong to the same part
of the partition $P=P_{I_\odd}\cup P_{I_\even}$ precisely when $[i,i']\cap J$ is of odd size. This is in  contradiction
with the criterion for $H_J$ to separate $p_i$ from $p_{i'}$. 
\end{proof}

By Ramsey's theorem, there is a number $\OT_d(n)$ such that each sequence of $\OT_d(n)$ points in general position contains 
an orientation-homogeneous subsequence of length~$n$. The growth rate of $\OT_d(n)$ is known quite precisely: For all $d\ge 2$ we have $\tw_d(c'_d n) \le \OT_d(n) \le \tw_d(c_d n)$ for positive constants $c'_d < c_d$. The upper bound is due to Suk~\cite{suk_ramsey}, and the lower bound is due to B\'ar\'any, Matou{\v{s}}ek and P\'or \cite{bmp}, which is based on an earlier work by Eli{\'a}{\v{s}}, Matou{\v{s}}ek, Rold{\'a}n-Pensado and Safernov{\'a} \cite{emrs}.

\subsection{Point selection}\label{subsec_pointselect}

The following lemma is a minor variation on Lemma~2.2 from \cite{aknss}:

\begin{lemma}\label{lem:tverberg_conv}
Let $s\eqdef \lfloor d/2\rfloor +1$, and let $D\eqdef (s-1)(d+1)+1$. Let $(p_1, p_2, \ldots, p_{2D+1})$ be an orientation-homogeneous sequence of $2D+1$ points in $\R^d$. Let $Q=\{p_2, p_4, \ldots, p_{2D}\}$ and $R = \{p_1, p_3, \ldots, p_{2D+1}\}$ (so $\abs{Q}=D$ and $\abs{R}=D+1$). Then $\Tver_s(Q)\in \conv R$.
\end{lemma}
\begin{proof}
Let $x\eqdef \Tver_s(Q)$. If $x\notin \conv R$, then there exists a hyperplane $H$ separating $x$ from $R$.
There must be at least $s$ points of $Q$ on the same side of $H$ as $x$ (at least
one from each part in the Tverberg partition). Let $Q'$ be any $s$ of these points. 
Pick any set $R'\subset R$ of $\lceil d/2\rceil +1$ points that interleaves $Q'$.
By Lemma~\ref{lem_interleave_radon}, the sets $\conv Q'$ and $\conv R'$ intersect, contradicting the fact that that $H$ separates $Q'$ from~$R$.
\end{proof}

\subsection{A regularity lemma for semialgebraic predicates}\label{sec_regularity}

We shall use a regularity lemma of Fox--Pach--Suk \cite{fox_pach_suk_regularity}, which is a quantitative improvement
over the prior version due to Fox--Gromov--Lafforgue--Naor--Pach \cite{FGLNP}. The improvement is due to the use of the efficient cuttings 
of Chazelle--Friedman \cite{chazelle_friedman} and Clarkson \cite{clarkson}.

Consider a polynomial $f\in \R[\vec{x}_1,\dotsc,\vec{x}_k]$, where each $\vec{x}_i$ is
a vector of $d$ indeterminates. The \emph{degree} of $f$ in $\vec{x}_i$ is the degree
of $f$ as a polynomial in $\vec{x}_i$ while regarding $\vec{x}_j$ for $i\neq j$ as constants.
We say that $f$ is of \emph{complexity} at most $D$ if it 
is of degree at most $D$ in each of $\vec{x}_1,\dotsc,\vec{x}_k$. 

\begin{lemma}[Theorem 1.3 in \cite{fox_pach_suk_regularity}]\label{lem:regularity}
For any $k,d,t,D\in \N$ there exists a constant $c=c(k,d,t,D)>0$
with the following property. Let $0<\gamma<1/2$, let $P\subset \R^d$ be a finite multiset,
and let $f_1,\dotsc,f_t\in \R[\vec{x}_1,\dotsc,\vec{x}_k]$ be $t$ polynomials of complexity at most
$D$ each. Then there exists a partition $P=P_1\cup\dotsb\cup P_M$ of $P$ into at most
$M\leq (1/\gamma)^{c}$ parts, and a small set $\mathcal{E}\subset [M]^k$ of ``exceptional" $k$-tuples, satisfying the following:
\begin{enumerate}
\item The exceptions are few: $\abs{\mathcal{E}}\leq \gamma M^k$,
\item\label{pt:regular} Almost all $k$-tuples are regular: whenever $(i_1,\dotsc,i_k)\notin \mathcal{E}$ and
$p_1\in P_{i_1},\dotsc,p_k\in P_{i_k}$, then the sign of
\[
  f_j(p_1,\dotsc,p_k)
\]
depends only on $j$ and on the tuple $(i_1,\dotsc,i_k)$ but not on the actual choice of the points $p_1,\dotsc,p_k$. (Note that
the elements $i_1,\dotsc,i_k$ of the tuple need not be distinct nor in increasing order.)
\item The partition is an equipartition: For all $i,j$ the cardinalities of $P_i$ and of $P_j$ differ by at most one.
\end{enumerate}
\end{lemma}

(The statement appearing in \cite{fox_pach_suk_regularity} is slightly different: In part~(\ref{pt:regular}) 
instead of claiming that the signs of all $f_j$ are constant, the original merely states that an arbitrary fixed Boolean formula in
signs of $f_j$ is constant. However, their proof actually establishes the stronger statement above. Alternatively, one may
refine the partition $\mathcal{P}$ by iterative application of the original statement to each $f_j$ in turn.
The only minor drawback is that instead of a true equipartition one would then obtain a partition
whose parts differ by as much as $t$, the number of polynomials.)

The main point of Lemma~\ref{lem:regularity} is that the number $M$ of parts is independent of $\abs{P}$
(otherwise we could trivially partition $P$ into parts of size $1$). The price for this independence is the small set $\mathcal{E}$ which indexes ``irregular'' tuples.

Invoking Lemma~\ref{lem:regularity} with the orientation predicate, we obtain the following result, which is what we actually need:
\begin{corollary}\label{cor:regularity}
For each $d$ there exists a constant $c=c(d)>0$ with the following property. Let $0<\gamma<1/2$, and let $P\subset \R^d$ be a finite point set in general position.
Then there exists a partition $P=P_1\cup\dotsb\cup P_M$ of $P$ into
$M$ parts, with $1/\gamma \leq M \leq 2(1/\gamma)^{c}$, and a small hypergraph $\Hg\subset \binom{[M]}{d+1}$ of ``exceptional" $(d+1)$-sets, satisfying the following:
\begin{enumerate}
\item $\abs{\Hg}\leq \gamma \binom{M}{d+1}$,
\item Whenever $\{i_0,i_1,\dotsc,i_d\}\in\binom{[M]}{d+1}\setminus \Hg$ and
$p_0\in P_{i_0},p_1\in P_{i_1},\dotsc,p_d\in P_{i_d}$, then the sign of
$\orient(p_0,p_1,\dotsc,p_d)$
depends only on the tuple $(i_0,\dotsc,i_d)$ but not on the actual choice of the points $p_0,\dotsc,p_d$. (The sign of $\orient$ obviously does depend
on the permutation of the elements $i_0,\dotsc,i_d$.)
\item For all $i,j$ the cardinalities of $P_i$ and of $P_j$ differ by at most one.
\end{enumerate}
\end{corollary}
\begin{proof}
If $\abs{P}\leq 2(1/\gamma)^c$, then simply partition $P$ into parts of size~$1$. So, assume $\abs{P}\geq 2(1/\gamma)^c$. We apply
Lemma~\ref{lem:regularity} with $t\eqdef 1$ and $f_1\eqdef \orient$. We obtain a partition of $P$ into $M\le (1/\gamma)^c$ parts, each of size at least~$2$, and a set $\mathcal{E}\subset [M]^{d+1}$ of size at most $\gamma M^{d+1}$.

We now show that all tuples $(i_0,\dotsc,i_d)\in [M]^{d+1}$ that contain repeated elements must belong to $\mathcal{E}$. Indeed, consider one such tuple, and say $i_j=i_{j'}$. Since the part $P_{i_j}$ has at least two elements, say $p$ and $q$, swapping $p$ and $q$ causes $\orient$ to flip its nonzero sign (recall that $P$ is in general position). Hence, the tuple $(i_0,\dotsc,i_d)$ is not regular, i.e., it does not satisfy property~2 above.

This consideration implies the lower bound for $M$: We have $\gamma M^{d+1} \ge \abs{\mathcal{E}}\geq M^{d+1}-(d+1)!\binom{M}{d+1}\geq M^{d+1}-M^d(M-1) = M^d$,
and hence $M\geq 1/\gamma$.

Finally, we let $\Hg$ consist of all tuples in $\mathcal{E}$ whose elements are pairwise distinct (this definition makes sense since, in our case, $\mathcal{E}$ is invariant under permutations). Since $\mathcal{E}$ contains \emph{all} the tuples with repeated elements, it can contain at most a $\gamma$-fraction of the remaining tuples. Therefore, the same is true for $\Hg$.
\end{proof}

\section{Independent sets in hypergraphs}\label{sec_hypergraph}

We will also need the following bound on hypergraph Tur\'an numbers. We give a simple probabilistic proof based on \cite[Theorem 3.2.1]{alon_spencer},
though a stronger bound can be found in \cite{decaen_turan}.

\begin{lemma}\label{lem_turan}
Let $r\geq 2$, and suppose $\Hg$ is an $r$-uniform hypergraph on~$n$ vertices with $\beta n^r$
edges, where $n\geq \tfrac{1}{2}\beta^{-1/(r-1)}$. Then $\Hg$ contains an independent set on at least $\tfrac{1}{4}\beta^{-1/(r-1)}$ vertices. 
\end{lemma}

\begin{proof}
Let $p\eqdef\beta^{-1/(r-1)}/(2n)$. Note that $p\leq 1$ by the assumption on~$n$.
Let $S\subseteq V(\Hg)$ be a random set where $\Pr[v\in S]=p$ for each $v\in V(\Hg)$
independently. Then the expected number of edges spanned by $S$ is $p^r\beta n^r$. For each edge
in $S$ we may remove one vertex to obtain an independent set. Hence,
$\Hg$ contains an independent set of size at least $\mathbb{E}[I]=pn-p^r \beta n^r\geq \tfrac{1}{2}\beta^{-1/(r-1)}-\tfrac{1}{2^r}\beta^{1-r/(r-1)}/n^r\cdot n^r$.
\end{proof}

\section{Interval chains}\label{sec_interval}
We will reduce the geometric problem of constructing one-sided \vepsapproximants/ to a combinatorial problem about interval
chains. Let $[i,j]$ denote the interval of integers $\{i,i+1,\dotsc,j\}$. We still write $[t]$ for $\{1,2,\dotsc,t\}$. 
An \emph{interval chain} of size $k$ (also called $k$-chain) in $[t]$ is a sequence of $k$ consecutive, disjoint, nonempty intervals
\begin{align*}
  I&\eqdef [a_1,a_2-1][a_2,a_3-1]\dotsb[a_k,a_{k+1}-1],
\end{align*}
where $1\le a_1<a_2<\dotsb<a_{k+1}\le t+1$. Interval chains were introduced by Condon and Saks \cite{condon_saks}. They were subsequently
used by Alon, Kaplan, Nivasch, Sharir and Smorodinsky \cite{aknss} and by Bukh, Matou\v{s}ek, Nivasch \cite{bukh_matousek_nivasch_epsilonnets} to obtain bounds for weak \vepsnets/ for orientation-homogeneous point sets. 

A $D$-tuple of integers $(x_1, \ldots x_D)$ is said to \emph{stab} a $k$-chain $I$ if each $x_i$ lies in a different interval of $I$.

The problem considered in~\cite{aknss} was to build, for given $D$, $k$, and $t$, a small-sized family $\mathcal F$ of $D$-tuples that stab all $k$-chains in $[t]$. Phrased differently, for each interval chain $I$ with \emph{at least} $k$ intervals, there should be \emph{at least} one $D$-tuple in $\mathcal F$ that stabs $I$.

In contrast, here we will consider the following problem: Given $D$, $\veps$, and $t$, we want to build a small-sized family (multiset) $\mathcal F$ of $D$-tuples such that, for each interval chain $I$ in $[t]$, if $\alpha t$ is the number of intervals in $I$, then at least an $(\alpha-\veps)$-fraction of the $D$-tuples in $\mathcal F$ stab $I$. We call such an $\mathcal F$ an \emph{$\veps$-approximating family}.

Our construction of \vepsapproximating/ families is similar to the statement of the regularity lemma for words, due to Axenovich, Person and Puzynina \cite{words_regularity}.
The lemma, which was also independently discovered by Feige, Koren and Tennenholtz \cite{feige_repeat} under the name of `local repetition lemma', can be used directly to construct
\vepsapproximating/ families. Doing so yields a family whose size is exponential in $1/\veps$. In contrast, we avoid using the full strength of the regularity lemma
and obtain a construction of polynomial size.

\begin{lemma}\label{lem:interval_chains}
Suppose $D\ge 2$ and $0<\veps<1$. Let $K\eqdef \lceil (D-1)\ln (4/\veps)\rceil$ and $t \eqdef m (D-1)^K$ for some integer $m\ge 4/\veps$. Then there exists an $\veps$-approximating family $\mathcal F$ of $D$-tuples in $[t]$, of size $\abs{\mathcal F} \le t$.
\end{lemma}

\begin{proof}
The argument is more conveniently phrased in the ``dual'' setting, in which $D$-tuples become $(D-1)$-interval chains and $\ell$-interval chains become $(\ell+1)$-tuples. 
Namely, the $D$-tuple $(x_1,\dotsc,x_D)$ becomes the interval chain $[x_1+1,x_2]\dotsb[x_{D-1}+1,x_D]$, and the
$\ell$-interval chain $[a_1,a_2-1]\dotsb[a_{\ell},a_{\ell+1}-1]$ becomes the tuple $(a_1,a_2,\dotsc,a_{\ell+1})$. Then a $(D-1)$-chain $C$ ``stabs'' a tuple $T$ if $T$ contains points on both sides of $C$, as well as inside each interval of $C$.

For each $k=0, 1, \ldots, K-1$, we partition $[t]$ into disjoint intervals of length $(D-1)^k$, by letting $B_{k,i} \eqdef \bigl[(i-1)(D-1)^k+1, i(D-1)^k\bigr]$ for $1\le i \le t/(D-1)^k$. Then we group these intervals into disjoint $(D-1)$-chains, by letting
\begin{equation*}
\mathcal F_k \eqdef \{B_{k,(i-1)(D-1)+1}\cdots B_{k,i(D-1)} : 1\le i\le t/(D-1)^{k+1} \}.
\end{equation*}
We call each $\mathcal F_k$ a \emph{layer}. Note that each chain in $\mathcal F_k$ fits exactly in an interval of layer $k+1$.

Then we define the multiset $\mathcal F$ by taking $w_k$ copies of $\mathcal F_k$ for each $0\le k\le K-1$, where
\begin{equation*}
w_k \eqdef (D-2)^k.
\end{equation*}
Hence, letting $E\eqdef (D-2)/(D-1)$,
we have $\abs{\mathcal F} = \sum_{k=0}^{K-1} w_k\abs{\mathcal F_k} = t(1 - E^K)$. Therefore, by the choice of $K$,
\begin{equation*}
t/2\leq \abs{\mathcal F} < t
\end{equation*}

Let $J$ be a subset of $[t]$, and let $\alpha t$ be the size of $J$. We claim that at least an $(\alpha-\veps)$-fraction of the chains in $\mathcal F$ stab $J$.

Call a $(D-1)$-chain $C\in\mathcal F$ \emph{empty} if $J$ does not intersect any interval of $C$, and \emph{occupied} otherwise. If $C$ is occupied, then call it \emph{fully occupied} if $J$ intersects all intervals of $C$, and \emph{partially occupied} otherwise.

For each $0\le k\le K-1$, let $\beta_k$ denote the fraction of chains of $\mathcal F_k$ that are occupied by $J$, and let $\gamma_k \le \beta_k$ denote the fraction of chains of $\mathcal F_k$ that are partially occupied by $J$.

\begin{claim}\label{claim_gamma}
For each $k$ we have $\beta_k \ge \alpha + (\gamma_0 + \cdots + \gamma_k)/(D-1)$.
\end{claim}

\begin{proof}
For each layer $j$, Let $\mathcal F'_j$ be the set of occupied chains of $\mathcal F_j$, and let $\mathcal F''_j\subset \mathcal F'_j$ be the set of those that are only partially occupied. Hence, $\mathcal F''_j$ covers a $\gamma_j$-fraction of $[t]$. From each chain $C\in\mathcal F''_j$ choose an empty interval, and let $\mathcal B_j$ be the union of these empty intervals. Hence, $\mathcal B_j$ covers a $(\gamma_j/(D-1))$-fraction of $[t]$. Furthermore, since each chain in $\mathcal F''_j$ contains a point of $J$, the sets $\mathcal B_0, \ldots, \mathcal B_k$ must be pairwise disjoint, as well as disjoint from $J$, and their union $\mathcal U \eqdef \mathcal B_0 \cup \cdots \mathcal B_k \cup J$ must be completely contained in the union of $\mathcal F'_k$. Hence, $\mathcal F'_k$ covers at least an $\bigl( \alpha + (\gamma_0 + \cdots + \gamma_k)/(D-1)\bigr)$-fraction of $[t]$, and the claim follows.
\end{proof}

Let us now derive a lower bound on the number of fully occupied chains in $\mathcal F$. By some tedious calculations we obtain:
\begin{align*}
\sum_{k=0}^{K-1} (\beta_k-\gamma_k) w_k\abs{\mathcal F_k} &\ge \sum_{k=0}^{K-1} \Bigl(\alpha + \frac{\gamma_0 + \cdots + \gamma_k}{D-1} -\gamma_k\Bigr)w_k\abs{\mathcal F_k}\\
&=\alpha\abs{\mathcal F}+\sum_{k=0}^{K-1}\gamma_k\left(\frac{1}{D-1}\sum_{j=k}^{K-1} w_j\abs{\mathcal F_j}-w_k\abs{\mathcal F_k}\right)\\
&=\alpha\abs{\mathcal F}-t \frac{E^K}{D-1} \sum_{k=0}^{K-1}\gamma_k\ge \alpha\abs{\mathcal F}-t E^K\beta_{K-1}\\
&\ge\left(\alpha-2E^K\right)\abs{\mathcal F}\ge (\alpha-\veps/2)\abs{\mathcal F};
\end{align*}
where the upper bound for $\sum\gamma_k$ was obtained from Claim~\ref{claim_gamma}.

Finally, note that in each layer $\mathcal F_k$ there are at most two fully occupied chains that do not stab $J$. Since $\abs{\mathcal F_k} \ge m \ge 4/\veps$, the said chains constitute at most an $(\veps/2)$-fraction of $\mathcal F$.
\end{proof}

\section{Construction of the one-sided approximants}\label{sec_construction}
In this section we prove Theorem~\ref{thm_main}. 

Let $s$ and $D$ be as in Lemma~\ref{lem:tverberg_conv}. Then let $t$ be as small as possible to satisfy the condition of Lemma~\ref{lem:interval_chains} with $\veps/2$ in place of $\veps$ (so $t$ is polynomial in $1/\veps$). Then define
\begin{equation}\label{eq:construction_parameters}
u\eqdef \lceil 4/\veps\rceil,  \qquad n_0 \eqdef tu, \qquad N \eqdef\OT_d(n_0),\qquad
\beta\eqdef (4N)^{-d},\qquad\gamma\eqdef \beta(\veps/5)^{d+1};
\end{equation}
where the function $\OT_d(n_0)$ is defined at the end of Section~\ref{subsec:orientation_homog}. Invoking Lemma~\ref{lem:interval_chains}, let $\mathcal{F}$ be an \vepstwoapproximating/ family of
$D$-tuples in $[t]$, of size $\abs{\mathcal{F}}\le t$.

Let $P\subset \R^d$ be a given finite point set, and let $n\eqdef\abs{P}$. We will construct a one-sided \vepsapproximant/ multiset $A$ for $P$.
If $n\leq 40/(\veps \gamma^c)$ for the constant $c$ of Corollary~\ref{cor:regularity}, then simply let $A\eqdef P$.
Hence, assume $n\geq 40/(\veps \gamma^c)$. In this case, our multiset $A$ will consist of Tverberg points of certain $D$-tuples of points of $P$.

We first handle the case when $P$ is in general position; then we handle degeneracies with a simple
perturbation argument. Hence, suppose the point set $P\subset \R^d$ is in general position (no $d+1$
points are affinely dependent).

We start by invoking Corollary~\ref{cor:regularity} on  $P$ and the parameter $\gamma$ given in \eqref{eq:construction_parameters}. 
We obtain a partition of $P$ into $1/\gamma\leq M\leq 2(1/\gamma)^c$ almost-equal-sized \emph{parts} $P_1,\dotsc,P_M$,
and a corresponding hypergraph $\Hg\subseteq \binom{[M]}{d+1}$ of size
$\abs{\Hg}\leq \gamma\binom{M}{d+1}$. 

We make all parts have exactly the same size by discarding at most one point from each part. Hence we discard at most
$M\leq 2(1/\gamma)^c$ points. Since $n\geq 40/(\veps \gamma^c)$, we discarded at most an $(\veps/20)$-fraction of the points of~$P$.
By a slight abuse of notation, we denote the new parts by the same names $P_1,\dotsc,P_M$.
We will consider $P_1,\dotsc,P_M$ as an \emph{ordered} sequence (where the order was chosen arbitrarily).

Let $\widehat P = (p_1, \ldots, p_M)$, where $p_i\in P_i$ for all $i$, be an arbitrarily chosen sequence of representatives from the parts. We will now repeatedly ``fish out'' equal-length orientation-homogeneous subsequences from $\widehat P$, until there are too few points left to continue the process. For this purpose, let $\widehat P_1 \eqdef \widehat P$, and let $i\leftarrow 1$. Repeat the following: If $\bigl|\widehat P_i\bigr| < \veps M/5$ then stop. Otherwise, $\widehat P_i$ is large enough so that the number of edges of $\Hg$ spanned by $\widehat P_i$ is at most
\begin{equation*}
|\Hg| \leq \gamma \binom{M}{d+1} \leq \gamma M^{d+1} = \beta (\veps M/5)^{d+1} \leq \beta \bigl\lvert\widehat P_i\bigr\rvert^{d+1}.
\end{equation*}
In view of $M\geq 1/\gamma$, we also have $\veps M/5\geq (5/\veps)^d/\beta \ge \frac{1}{2}\beta^{-1/d}$.
Hence, we can apply Lemma~\ref{lem_turan} on $\widehat P_i$ with $r=d+1$. We conclude that $\widehat P_i$ has an independent set of size $N$. By the definition of $N$, that independent set has an orientation-homogeneous subsequence $S_i$ of length $n_0$. Let $\widehat P_{i+1} \eqdef \widehat P_i \setminus S_i$, increase $i$ by $1$, and return to the beginning of the loop.

At the end of this process, we get orientation-homogeneous sequences $S_1,\allowbreak S_2,\allowbreak \ldots,\allowbreak S_m$ for some $m\leq M/n_0$, and a leftover sequence $S^*\eqdef \widehat P_{m+1}$ of size at most $\veps M/5$. From each $S_i$ we will now construct a multiset~$A_i$ of Tverberg points; their union will be our desired multiset~$A$.

So fix $i$, and denote $S_i = (q_0,q_1, q_2, \ldots, q_{n_0-1})$. Let $v_j \eqdef q_{(j-1)u}$ for all $1\le j\le t$. We will call the elements $v_j$ \emph{separators}. Let $\mathfrak v\eqdef (v_1, \ldots, v_t)$.
For each $j=1,\dotsc,t$, define the \emph{block} $b_j \eqdef (q_{(j-1)u+1}, \ldots, q_{ju-1})$, which contains the elements of $S_i$ between separators $v_j$ and $v_{j+1}$. Let 
\[
  \mathcal{B}_j\eqdef\bigcup_{p_k\in b_j}  P_k
\]
be the union of all the parts that correspond to points of block~$b_j$.

To each $D$-tuple $\overline x = (x_1, \ldots, x_D)\in \mathcal F$, associate the $D$-tuple of separators $Q_{\overline x}\eqdef \{v_{x_1},\dotsc,v_{x_{D}}\}$. Then define the multiset
\[
  A_i\eqdef \{ \Tver_s(Q_{\overline x}) : \overline x\in \mathcal{F} \}.
\]

\begin{lemma}\label{lem:regul_use}
Let $C\subseteq \R^d$ be a convex set. Take the set of indices $J\eqdef \{j :  \mathcal{B}_j\cap C\neq\emptyset\}$. List the elements of $J$ in increasing order as $J=\{j_1, j_2, \ldots, j_\ell\}$. Let $I$ be the $(\ell-1)$-interval chain:
\begin{equation*}
I \eqdef [j_1+1, j_2][j_2+1,j_3]\cdots[j_{\ell-1}+1,j_\ell].
\end{equation*}
Then, if the $D$-tuple $\overline x\in\mathcal F$ stabs $I$, then $C$ contains the corresponding Tverberg point $\Tver_s(Q_{\overline{x}})$.
\end{lemma}

\begin{proof}
Suppose $\overline x = (x_1, \ldots, x_D)$ stabs $I$. Then there exists a subset $J'\eqdef (j'_0, \ldots, j'_D)\subset J$ such that $j'_0 < x_1 \le j'_q < \cdots < x_D \le j'_D$. For each $j \in J$ there is a part $P(j)$ whose representative point $p(j)$ belongs to the block $b_j$, and such that $C$ contains some point $p'(j)\in P(j)$. The sequence of representatives $p(j'_0), v_{x_1}, p(j'_1), \ldots, v_{x_D}, p(j'_D)$, being a subsequence of $S_i$, is orientation-homogeneous. Therefore, by regularity, and since $S_i$ avoids the hypergraph $\Hg$, the sequence  $p'(j'_0), v_{x_1}, p'(j'_1), \ldots,\allowbreak v_{x_D}, p'(j'_D)$ is also orientation-homogeneous. Therefore, by Lemma~\ref{lem:tverberg_conv}, we have
\begin{equation*}
\Tver_s(Q_{\overline{x}}) \in \conv{\{p'(j'_0),\allowbreak \ldots, p'(j'_D)\}}\subset C,
\end{equation*}
as desired.
\end{proof}

Let $\mathcal{S}_i\eqdef \bigcup_{p_j\in S_i} P_j$ be the union of all the parts whose representative points belong to $S_i$.

\begin{corollary}\label{cor_si}
Let $C\subseteq \R^d$ be a convex set, and let $\alpha$ be the fraction of the points of $\mathcal{S}_i$ contained in~$C$. Then $C$ contains at least an
$(\alpha-3\veps/4)$-fraction of the points of~$A_i$.
\end{corollary}
\begin{proof}
Since $\abs{\mathfrak{v}}=t\leq \veps n_0/4$, and since all the parts $P_1,\dotsc,P_M$ have equal size,
the set $C$ meets at least an $(\alpha-\veps/4)$-fraction of the sets $\mathcal{B}_j$. The desired conclusion follows 
from Lemma~\ref{lem:regul_use} since $\mathcal{F}$ is \vepstwoapproximating/.
\end{proof}

Finally, let
\[
  A\eqdef \bigcup_{i=1}^m A_i.
\]
With some patience, we can use \eqref{eq:construction_parameters} and the bound $\OT_d(n)\leq \tw_d(c_d n)$ mentioned above
to obtain the bound $\abs{A}=M\abs{\mathcal{F}}\leq \tw_d\bigl(\veps^{-c'}\bigr)$ for some constant $c'=c'(d)>1$.

Note that at most $(\veps/20)n+(\veps/5)n=\veps n/4$
points of $P$ were either discarded in making $P_1,\dotsc,P_M$ equal or were relegated to the ``leftover''
$S^*=\widehat{P}\setminus (S_1\cup \dotsb \cup S_m)$. So, if a convex set $C$
contains an $\alpha$-fraction of the points of $P$, and an $\alpha_i$-fraction of the points of $\mathcal{S}_i$
for each $i$, then $\avg_i \alpha_i \geq \alpha-\veps/4$. 

By Corollary~\ref{cor_si}, $C$ contains at least an $(\alpha_i-3\veps/4)$-fraction of the points of~$A_i$.
Hence, averaging again, $C$ contains an $(\alpha-\veps)$-fraction of the points of~$A$.

This concludes the proof of Theorem~\ref{thm_main} for the case when $P$ is in general position.\medskip

If $P = \{p_1, \ldots, p_n\}$ is not in general position, take an arbitrarily small continuous perturbation $P(t) \eqdef \{p_1(t), \ldots, p_n(t)\}$ such that $P(0) = P$ and $P(t)$ is in general position for all $0<t\le 1$. For each $t>0$ we apply the above argument on $P(t)$; we get a family $\mathcal I(t)\subset \binom{[n]}{D+1}$ such that multiset $A(t) \eqdef \{\Tver_s(P(t)_I) : I \in \mathcal I(t)\}$ is a one-sided \vepsapproximant/ for $P(t)$. Since $P$ is finite, there are only a finitely many possible values for $\mathcal I(t)$, so one of them occurs infinitely often for $t = t_1, t_2, t_3, \ldots$ with $\lim t_i = 0$. Then, by a standard argument, the limit multiset $\lim_{i\to\infty} A(t_i)$ exists and is a one-sided \vepsapproximant/ for $P$.

\section{Problems and remarks}\label{sec_remarks}

\begin{itemize}
\item The main problem is to prove reasonable upper bounds on $g(\veps,d)$. The only known lower
bound on $g(\veps,d)$ is of the form $c_d(1/\veps)\log^{d-1}(1/\veps)$.
It is a consequence of the lower bounds on the size of weak \vepsnets/ \cite{bukh_matousek_nivasch_epsilonnets} and
the fact that every one-sided \vepsapproximant/ is an \vepsnet/.

\item Much smaller one-sided approximants can be constructed if $P$ is orientation-homogeneous: We apply the same construction that was applied to individual sets $S_i$ in Section~\ref{sec_construction} to the set $P$
(with $u\eqdef \abs{P}/t$ instead of $u=n_0/t$), obtaining
one-sided \vepsapproximants/ of size polynomial in $1/\veps$. While this bound is much better than the general bound on $g(\veps,d)$ from
Theorem~\ref{thm_main}, it is still far from the known bounds for \vepsnets/: 
Every orientation-homogeneous set admits an \vepsnet/ of size only $O\bigl(\veps^{-1}\alpha(\veps^{-1})\bigr)$ in the plane
and of size only $\veps^{-1}2^{\alpha(\veps^{-1})^{O(1)}}$ in $\R^d$ for $d\geq 3$,
where $\alpha$ is the inverse Ackermann function \cite{aknss}.

\item The \emph{diagonal of the stretched grid} is a specific orientation-homogeneous sequence considered in \cite{bukh_matousek_nivasch_stabbing} and in~\cite{bukh_matousek_nivasch_epsilonnets}. Denote it $D$. The authors in~\cite{bukh_matousek_nivasch_epsilonnets} obtained a lower bound for \vepsnets/ for $D$ from the lower bound for the interval chains problem considered in \cite{aknss}. Similarly, a lower bound for the interval chains problem discussed in Section~\ref{sec_interval} would yield a lower bound for \vepsapproximants/ for $D$.

\item In Theorem~\ref{thm_main} it is possible to assure that the one-sided approximant~$A$ is a genuine set rather than a multiset.
It is easy to do so if $P$ is in general position, as we may simply perturb each point of $A$ slightly. In general,
we cannot ensure that each sequence $S_i$ is orientation-homogeneous, but we can ensure that
each $S_i$ is orientation-homogeneous inside the affine subspace $\ahull S_i$. That can be done by
using Ramsey's theorem to extract subsequences of $\widehat{P}$ that lie in a proper affine subspace,
and then using the induction on the dimension. We can then perturb the points of $A_i$ within $\ahull S_i$. The rest of the argument remains the same.
\end{itemize}

\textbf{Acknowledgement.} We thank the referees for pointing omissions and typos. We are also grateful to
Po-Shen Loh and the late Jirka Matou\v{s}ek for discussions.

\bibliographystyle{alpha}
\bibliography{onesidedapprox}

\end{document}